# ON THE ASYMPTOTIC EXPANSIONS FOR TIME-VARYING SCALAR DIFFERENTIAL EQUATIONS POSSESSSING LIMITING DIFFERENTIAL EQUATIONS BY APPLICATION OF THE RESIDUE THEOREM TO THEIR DISCRETIZED COUNTERPARTS. PRELIMINARY RESULTS

M. De la Sen
Department of Electricity and Electronics. Faculty of Science and Technology
Campus of Leioa (Bizkaia)- SPAIN

## 1. Introduction

Asymptotic expansions for higher-order scalar difference equations about their limiting equations are obtained in [1] via the residue theorem and the z-transform. The results derived in [1] are nearly dual to previous ones obtained for ordinary continuous-time linear time-varying and for functional differential equations. See [4] and references therein. In this paper, it is proved that, in a dual context, asymptotic expansions of ordinary linear time-differential equations which possess limiting equations to their limiting equations might be obtained by first discretizing them and then using residues calculation for such discrete- time counterparts according to the formulation provided in [1]. The constructive derivation of the necessary technical proofs is performed by using the residues theorem on the discrete-time associated difference equations. For such a purpose, a time-dependent time-varying finite sampling period is defined for each time instant at which the asymptotic expansion of the solution of the continuous-time differential equation is investigated. The (time-varying) bounded sampling period is defined as the quotient of each current time instant and a positive integer which tends to infinity as time tends to infinity. Such a sampling period is only used for obtaining the residue formulas allowing to build the continuous-time solution for each time instant from the discrete-time difference equation. The asymptotic expansion formulas may be directly extended to the case of stable Lyapunov´s equations and to the use of a constant sampling period. This second extension is possible since the errors of the solution trajectories with respect to their previously sampled values remain uniformly bounded for both the current differential equation and its limiting counterpart within the inter-sample period. Complementary stability results are also obtained through the manuscript.

**1.1 Notation**. $\boldsymbol{R}$ is the set of real numbers, $\boldsymbol{R}_+ := \{ z \in \boldsymbol{R} : z \geq 0 \}$, $\boldsymbol{R}_{0+} := \boldsymbol{R}_+ \cup \{ 0 \} = \{ z \in \boldsymbol{R} : z \geq 0 \}$.

$\boldsymbol{C}$ is the set of complex numbers, $\boldsymbol{C}_+ := \{ z \in \boldsymbol{C} : Re\, z \geq 0 \}$, $\boldsymbol{C}_{0+} := \boldsymbol{R}_+ \cup \{ 0 \} = \{ z \in \boldsymbol{C} : Re\, z \geq 0 \}$.

**N** is the set of natural numbers , $\boldsymbol{N}_0 := \boldsymbol{N} \cup \{ 0 \}$ and $\overline{n} := \{ 1, 2, ..., n \}$; $\forall n \in \boldsymbol{N}$ .

Let $x(t) = \sum_{k=0}^{\infty} x(kT) \delta(t - kT)$ bea real-valued function or vector function of support of zero measure . $S_T := \{ r \in \boldsymbol{R}_{0+} : r = kT , \forall k \in \boldsymbol{N}_0 \}$ with $\delta(t - kT)$ being a Dirac impulse at $t = kT$ and $x_k \equiv x(kT) = x(t)\big|_{t=kT}$ and whose range is the real sequence $\{ x(kT) \}_0^{\infty}$. Then, $\hat{x}_T(z) := \sum_{k=0}^{\infty} x(kT) z^{-k}$ is the Z-transform of $x(t)$ defined over some definition domain in $\boldsymbol{C}$ provided that such a domain exits and it is not empty. Equivalently, it can be said that $\hat{x}_T(z)$ is the – transform of the real sequence , or vector sequence, $\{x(kT)\}_0^{\infty}$.

The functions $id : \boldsymbol{N}_0 \times \boldsymbol{C} \to \boldsymbol{C}$ and $id^{n-1} : \boldsymbol{C} \to \boldsymbol{C}$ are defined by $id(n-1, z) = id^{n-1}(z) = z^{n-1}$.

## 2. The continuous- time differential equation, its discrete-time counterpart and their limiting versions

Consider the nth-order Poincaré differential equation on $\boldsymbol{R}_{0+}$

$$\sum_{i=0}^{n} \alpha_{n-i}(t) D^i x(t) = 0\,,\ D^i x(0) = x_{0i} \in \boldsymbol{C} \quad (\forall i \in \overline{n-1} \cup \{0\}) \tag{2.1}$$

where $\bar{j} := \{1, 2, \dots, j\} \subset \boldsymbol{N}$ and $p(D,t) := \sum_{i=0}^{n} \alpha_{n-i}(t) D^i$ is a monic polynomial complex–valued function of domain $\boldsymbol{R}_{0+}$ in the indeterminate $D := d/dt$ which is the time-derivative operator subject to the recursive formulas $D^{i+1} = D.D^i$ with $D^0 = 1$, point-wise defined by $\left(D^i f\right)(t) = f^{(i)}(t) = \frac{d^i f(t)}{dt^i}$ ; $\forall f \in C^{(i-1)}\left(\boldsymbol{R}_{0+}, \boldsymbol{C}^q\right)$, for any given integer $q \ge 1$, such that $f^{(i)}(t)$ exists everywhere in $\boldsymbol{R}_{0+}$; and $\alpha_i : \boldsymbol{R}_{0+} \to \boldsymbol{C}$ ; $\forall i \in \bar{n} \cup \{0\}$ , with $\bar{n} := \{1, 2, \dots, n\}$, are piecewise continuous bounded functions with $\alpha_0(t) = 1$; $\forall t \in \boldsymbol{R}_{0+}$ (since $p(D,t)$ is monic). It is a direct to see by construction that (2.1) is equivalent to the time-varying – linear differential system:

$$\dot{z}(t) := D z(t) = A(t) z(t) \tag{2.2}$$

where $z(t) := \left(x(t), Dx(t), \dots, D^{n-1} x(t)\right)^T$, subject to the initial condition vector $z(0) = z_0 := \left(x_0, Dx(0), \dots, D^{n-1} x(0)\right)^T$, and

$$A(t) = \left(a_{ij}(t)\right) := \begin{bmatrix} 0 & \vdots & I_{n-1} \\ & & \\ -\alpha_n(t) & -\alpha_{n-1}(t) & \cdots\ -\alpha_1(t) \end{bmatrix} \tag{2.3}$$

since the functions $\alpha_i : \boldsymbol{R}_{0+} \to \boldsymbol{C}$ are piecewise continuous, it exists a unique solution on $\boldsymbol{R}_{0+}$ to the differential system (2.2), equivalently to the differential equation (2.1), for each $z(0) = \left(x_0, Dx(0), \dots, D^{n-1} x(0)\right)^T$ from Picard- Lindeloff theorem defined by

$$z(t) = \Psi(t, 0) z(0) \tag{2.4}$$

where $\Psi(t, \tau)$ is the unique evolution operator associated with (2.2) which satisfies:

$$\dot{\Psi}(t,\tau) := \frac{d\Psi(t,\tau)}{dt} = A(t)\Psi(t,\tau) = \Psi(t,\tau)A(t)\ ;\ \forall t, \tau \in \boldsymbol{R}_{0+} \tag{2.5}$$

subject to $\Psi(t,t)=I_n$; $\forall t\in\mathbf{R}_{0+}$; $\forall\tau,\forall t\in\mathbf{R}_{0+}$. Note that $\Psi(t,\tau)$ is not, in general, a convolution operator if $A(t)$ varies with time. If $A(t)$ is a constant matrix A, $\forall t\in\mathbf{R}_{0+}$ then $\Psi(t,\tau)=\Psi(t-\tau)=e^{A(t-\tau)}$ is a convolution operator and also a $C_0$- semigroup of infinitesimal generator A. A particular case of both theoretical and practical wide interest concerning (2.5) is addressed in the subsequent result.

THEOREM 2.1. The following properties hold:

(i) Assume that and that all the entries of $A(t)$ are in $L_1(\mathbf{R}_{0+},\mathbf{C})$ and that $A(t)\,e^{\int_0^t A(\tau)\,d\tau}$ commute. Then, $\Psi(t,t_1)=e^{\int_{t_1}^t A(\tau)\,d\tau}$ is the unique evolution operator of the differential system which satisfies $\Psi(t,0)=\Psi(t,t_1)\Psi(t_1,0)$ $\forall t\,,\forall t_1\in\mathbf{R}_{0+}$ and $\Psi(t_1,t)=\Psi^{-1}(t,t_1)$.

(ii) $\Psi((k+n)T,kT)=e^{\overline{A}_k nT}=\left(e^{\overline{A}_k T}\right)^n$ for any $T\in\mathbf{R}_+$ some constant $\overline{A}\in\mathbf{C}^{n\times n}$ dependent on the interval $[kT,(k+1)T)$ . Also, if $\lim\limits_{t\to\infty} A(t)=A^*=\left(a_{ij}^*\right)$ then

$$\lim_{k\to\infty} a_{ij}\left(\xi_{ij}(k)\right)=a_{ij}^*\,;\forall i,j\in\overline{n}\Rightarrow\lim_{k\to\infty} A_{\xi}(k)=A^*\Rightarrow\lim_{k\to\infty}\Psi(t,kT)=\lim_{k\to\infty} e^{A^*t}$$

$\forall t\in[kT,(k+1)T)$, which is a $C_0$- semigroup of infinitesimal generator $A^*$ being independent on each particular T.

(iii) The solution of the Poincaré differential system (2.4) is uniquely described at its samples associated with any given sampling period $T\in\mathbf{R}_+$ by a discrete- time difference system $w_{k+1}=\Phi_k w_k$; $\forall k(\geq n)\in\mathbf{N}$ for a vector of initial conditions $w_{n-1}=\left(x_{n-1},x_{n-2},\dots,x_0\right)^T=M_k z(0)$ and any vector $z(0)=\left(x_{00},\dots,x_{0,n-1}\right)^T$ of initial conditions of (2.4) , taken from those of the differential equation (2.1), with

$$\Phi_k:=M_{k+1}\Psi((k+1)T,kT)M_k^{-1}\;;\quad M_k:=E\overline{M}_k \tag{2.6a}$$

$$E:=\begin{bmatrix} e_1^T & 0^T & \cdots & 0^T\\ 0^T & e_1^T & \cdots & 0^T\\ 0^T & \cdots & \ddots & 0^T\\ 0^T & 0^T & \cdots & e_1^T\end{bmatrix}\in\mathbf{R}^{n\times n^2}\;;\;\overline{M}_k:=\begin{bmatrix}\Psi(kT,(k-n+1)T)\\ \Psi((k-1)T,(k-n+1)T)\\ \vdots\\ \Psi((k-n+2)T,(k-n+1)T)\\ I_n\end{bmatrix}\in\mathbf{R}^{n^2\times n} \tag{2.6b}$$

where $M_k\in\mathbf{R}^{n\times n}$ is non-singular; $\forall k(\geq n)\in\mathbf{N}$, $e_1=(1,0,\cdots,0)^T\in\mathbf{R}^n$ is the nth –unity Euclidean vector. If $\lim\limits_{t\to\infty} A(t)=A^*=\left(a_{ij}^*\right)$ then, for any given $T\in\mathbf{R}_+$, $w_{k+1}=\Phi_k w_k$ has a limiting difference system $w_{k+1}^*=\Phi_k^* w_k^*$ as $k\to\infty$ depending on each particular T.

*Proof*: (i) $\Psi(t, t_1)$ as defined is unique $\forall t, \forall t_1 \in \mathbf{R}_{0+}$ by simple inspection. Also,

$$\Psi(t,0)=e^{\int_{t_1}^{t} A(\tau)d\tau+\int_0^{t_1} A(\tau)d\tau}=\left(e^{\int_{t_1}^{t} A(\tau)d\tau}\right)\left(e^{\int_0^{t_1} A(\tau)d\tau}\right)=\Psi(t,t_1)\Psi(t_1,0) \qquad (2.7)$$

$\forall t, \forall t_1 \in \mathbf{R}_{0+}$ since the second identity in (2.6) holds since $\int_0^t A(\tau)d\tau$ and $\int_0^{t_1} A(\tau)d\tau$ commute . Also, direct calculations yield:

$$\Psi(t_1,t)=e^{\int_t^{t_1} A(\tau)d\tau}=e^{-\int_{t_1}^{t} A(\tau)d\tau}=\Psi^{-1}(t,t_1);\ \forall t_1(>t),\ \forall t \in \mathbf{R}_{0+} \qquad (2.8)$$

and

$$\dot{\Psi}(t,t_1)=A(t)e^{\int_{t_1}^{t} A(\tau)d\tau}=A(t)\Psi(t,t_1)=\Psi(t,t_1)A(t) \qquad (2.9)$$

which verifies (2.7). Thus, $\Psi(t, t_1)$, as defined, is an evolution operator $\forall t, \forall t_1 \in \mathbf{R}_{0+}$. Property (i) has been proven.

(ii) Note from Theorem 2.1(i) that

$$\begin{aligned} z(t) &= \Psi(t,0)z(0) \\ &= \Psi(t,t_1)\Psi(t_1,0)z(0)=\left(e^{\int_{t_1}^{t} A(\tau)d\tau}\right)z(t_1)\ ;\ \forall t(\ge t_1), \forall t_1 \in \mathbf{R}_{0+} \end{aligned} \qquad (2.10)$$

A sequence of (n+1) consecutive samples with any sampling period $T \in \mathbf{R}_+$ of (2.10) obeys the equation:

$$z_{k+n} := z\left[(k+n)T\right] = \Psi((k+n)T, (k+n-1)T)z_{k+n-1} = \Psi((k+n)T, kT)z_k \qquad (2.11)$$

and $x_{k+n} = e_1^T z_{k+n}$ where

$$\Psi((k+n)T, kT) = e^{\int_{kT}^{(k+n)T} A(\tau)d\tau} = e^{\sum_{i=1}^{D_k+1} A_{\xi_i}(k)T_i(k)} = e^{A_\xi(k)nT} \qquad (2.12)$$

provided that $A : \mathbf{R}_{0+} \to \mathbf{C}^{n\times n}$ has piecewise continuous entries with bounded discontinuities at times $kT + \sum_{i=1}^{D_k+1} T_i(k) \in [kT, (k+n)T)$ for some $D_k \in \mathbf{N}_0 := \mathbf{N} \cup \{0\}$ and some $T_i(k) \in \mathbf{R}_+$ $(\forall i \in \overline{D}_k)$ which are dependent on the time interval $[kT, (k+n)T)$, where $T_{D_k+1} = nT - \sum_{i=1}^{D_k} T_i$. The following identities hold in (2.11):

$A_\xi(k) := \left(a_{ij}(\xi_{ij}(k))\right) = \sum_{\ell=1}^{D_k+1} \frac{T_\ell}{nT} A_{\xi_\ell}(k)$ where $A_{\xi_\ell}(k) := \left(a_{ij}(\xi_{ij\ell}(k))\right)$ for some

$\xi_{ij\ell}(k) \in \left( kT + \sum_{i=1}^{\ell} T_i(k),\ kT + \sum_{i=1}^{\ell+1} T_i(k) \right)$; $\forall i, j \in \overline{n}$, $\forall \ell \in \overline{D}_k := \{1,2,....,D_k\}$ if $D_k \neq 0$;

and

$A_\xi(k) := \left(a_{ij}\left(\xi_{ij}(k)\right)\right) = \frac{\int_{kT}^{(k+n)T} A(\tau)\, d\tau}{nT}$ if $D_k = 0$ ( i.e. $A : \boldsymbol{R}_{0+} \to \boldsymbol{C}^{n\times n}$ is everywhere continuous within the real interval $\left(kT,\ (k+n)T\ \right)$) such that the identities below hold from La Rolle mean theorem for integrals with everywhere continuous integrals:

$$\int_{kT}^{(k+n)T} A(\tau)\, d\tau = \sum_{i=1}^{D_k} A_{\xi_i}(k) T_i(k) = A_\xi(k)\, nT$$

applied for each entry $a_{ij}(\tau)$ of the matrix function $A : \boldsymbol{R}_{0+} \to \boldsymbol{C}^{n\times n}$ on each time interval $\left( kT + \sum_{i=1}^{j} T_i(k),\ kT + \sum_{i=1}^{j+1} T_i(k) \right)$ leading to a, in general entry-dependent, mean value $a_{ij}\left(\xi_{ij}(k)\right)$ within such an interval. Note that (2.11) is a $C_0$- semigroup of infinitesimal generator $A_\xi(k)$ of the differential system $\dot{z}(t) = A_\xi(k) z(t)$ on the finite interval $\left[kT,\ (k+n)T\right)$ with initial conditions $z_k := z(kT)$. On the other hand, $\Psi(t, kT) = e^{A_\xi(k)t}$ ; $\forall t \in \left[kT, (k+1)T\right)$ is a fundamental matrix of such a differential system and then non singular, $\forall t \in \left[nT, (k+n)T\right)$ which fulfils $\Psi((k+n)T, kT) = e^{A_\xi(k)nT} = \left(e^{A_\xi(k)T}\right)^n$. The first part of Property (ii) has been proven with $\overline{A}_k := A_\xi(k)$. Furthermore, if $\lim_{t\to\infty} A(t) = A^* = \left(a_{ij}^*\right)$ then:

$$\lim_{k\to\infty} a_{ij}\left(\xi_{ij}(k)\right) = a_{ij}^* ; \forall i, j \in \overline{n} \Rightarrow \lim_{k\to\infty} A_\xi(k) = A^* \Rightarrow \lim_{k\to\infty} \Psi(t, kT) = \lim_{k\to\infty} e^{A^* t} \;\; ; \;\; \forall t \in \left[kT, (k+1)T\right),$$

which is a $C_0$- semigroup of infinitesimal generator $A^*$ which is independent on each particular T. Property (ii) has been proven.

(iii) $M_k := E\overline{M}_k$ is nonsingular since E and $\overline{M}_k$ are both nonsingular by construction, E possesses non less columns than rows and $\overline{M}_k$ has non less rows than columns, $\forall T \in \boldsymbol{R}_+$. Furthermore, $w_{n-1} = M_k z(0)$ is unique for each given $z(0)$ since $M_k$ is unique from its definition from given values the evolution operator $\Psi(kT, jT)$. Then, the existence and uniqueness of the solution $w_{k+1} = \Phi_k w_k ; \forall k\, (\geq n) \in \boldsymbol{N}$ under initial conditions $w_{n-1} = M_k z(0)$ follows from the identities

$$w_{k+1} = M_{k+1} z_{k+1} = M_{k+1} \Psi((k+1)T, kT) M_k^{-1} w_k \tag{2.13}$$

The solution is identical to that of (2.1), equivalently to that of (2.4) at sampling instants $t_k = kT$ since it is unique for any given $T \in \boldsymbol{R}_+$. From (2.5)-(2.6), note that $\exists \lim_{k\to\infty} \Phi_k = \Phi^*$, $\lim_{k\to\infty} \overline{M}_k = \overline{M}^*$, $\lim_{k\to\infty} M_k = M^*$ if $\lim_{t\to\infty} A(t) = A^* = \left(a_{ij}^*\right)$ implying that a limiting difference system $w_{k+1}^* = \Phi_k^* w_k^*$ exists as f $k \to \infty$ or each given $T \in \boldsymbol{R}_+$. □

Note that the existence of the limiting system referred to in Theorem 2.1 (iii) implies that a limiting scalar difference equation exists, [1],[5]. A more general result than Theorem 2.1 is the following one:

THEOREM 2.2. Assume that $A(t)=\boldsymbol{A}+\tilde{A}(t)$ and that all the entries of $\tilde{A}(t)$ are in $L_1\left(\boldsymbol{R}_{0+},\boldsymbol{C}\right)$. Then, the following properties hold:

(i) $\Psi(t,0)=e^{\boldsymbol{A}t}\left(I_n+\int_0^t e^{-\boldsymbol{A}\tau}\tilde{A}(\tau)\Psi(\tau,0)d\tau\right)$ (2.14)

is the unique evolution operator of the differential system $\forall(t\geq 0)\in\boldsymbol{R}_{0+}$. Furthermore, if $\lim_{t\to\infty} A(t)=A^*$ then $\exists \lim_{t\to\infty}\Psi(t,0):=\Psi^*(t,0)=e^{A^*t}=e^{\boldsymbol{A}t}\left(I_n+\int_0^t e^{-\boldsymbol{A}\tau}\tilde{A}^*\Psi^*(\tau,0)d\tau\right)$

where $\tilde{A}^*:=A^*-\boldsymbol{A}$.

(ii) Theorem 2.1 (iii) holds under the given assumptions.

*Proof:* (i) If (2.14) is true then

$$z(t)=e^{\boldsymbol{A}t}\left(I_n+\int_0^t e^{-\boldsymbol{A}\tau}\tilde{A}(\tau)\Psi(\tau,0)d\tau\right)z(0) \qquad (2.15)$$

On the other hand, one gets by taking time-derivatives in (2.15) with respect to the argument t that

$$\frac{d\,\Psi(t,t_1)}{dt}=\boldsymbol{A}e^{\boldsymbol{A}t}\left(I_n+\int_0^t e^{-\boldsymbol{A}\tau}\tilde{A}(\tau)\Psi(\tau,0)d\tau\right)+\tilde{A}(t)\Psi(t,0) \qquad (2.16)$$

Thus,

$$\dot{z}(t)=\boldsymbol{A}e^{\boldsymbol{A}t}\left(I_n+\int_0^t e^{-\boldsymbol{A}\tau}\tilde{A}(\tau)\,\Psi(\tau,0)d\tau\right)z(0)+\tilde{A}(t)\Psi(t,0)z(0)$$
$$=Az(t)+\tilde{A}(t)z(t) \qquad (2.17)$$

Conversely, if (2.16) integrated for any time, one gets $z(t)=e^{\boldsymbol{A}t}\left(I_n+\int_0^t e^{-\boldsymbol{A}\tau}\tilde{A}(\tau)\Psi(\tau,0)d\tau\right)z(0)+C$ for any $C\in\boldsymbol{C}$ but $z(t)\,]_{t=0}=z(0)$ for t = 0 and $e^{At}\,\rfloor_{t=0}=I_n$ imply $C=0$. Thus, (2.15) is true. The second part of the Theorem is proven as follows. Let $\Psi^*(t,0):=e^{A^*t}$ the $C_0$ -semigroup and evolution operator of $\dot{z}^*(t)=A^*z^*(t)$, assumed under initial conditions $z^*(0)=z(0)=z_0$, so that $\dot{\Psi}^*(t,0):=A^*e^{A^*t}$. Direct calculations yield for any matrix norm:

$$\left\|\dot{\Psi}(t,0)-\dot{\Psi}^*(t,0)\right\|\leq\left\|\boldsymbol{A}\,\Psi(t,0)+\tilde{A}(t)\Psi^*(t,0)-\boldsymbol{A}\,\Psi^*(t,0)-\tilde{A}^*\Psi^*(t,0)\right\|$$
$$\leq\left\|\boldsymbol{A}\left(\Psi(t,0)-\Psi^*(t,0)\right)\right\|+\left\|\tilde{A}(t)\Psi(t,0)-\tilde{A}^*\Psi^*(t,0)\right\| \qquad (2.18)$$

The proof is made by complete induction by assuming that $\left\| \Psi\left(t+i\delta,0\right)-\Psi^{*}\left(t+i\delta,0\right)\right\| \le \varepsilon_1\left(t+i\delta\right)$ for some given $\varepsilon_1 \in \boldsymbol{R}_+$ any arbitrary $\delta \in \boldsymbol{R}_+$ and $\forall t \ge t_1\left(\varepsilon_1\right)$, $\forall i \in \bar{j} \setminus \{0\}$ and some given $j \in \boldsymbol{N}_0$. Now, $\forall \varepsilon_1, \varepsilon_2 \in \boldsymbol{R}_+$, $\exists$ finite sufficiently large $t_i\left(\varepsilon_i\right)$; i =1,2 and strictly non-negative monotonically decreasing real sequences $\left\{\varepsilon_i\left(t+j\delta\right)\right\}_{j \in \boldsymbol{N}_0}$ for the given arbitrary $\delta \in \boldsymbol{R}_+$ satisfying $\varepsilon_i\left(t\right)=\varepsilon_i$; i =1,2 for any $t \ge max\left(t_1\left(\varepsilon_1\right), t_2\left(\varepsilon_2\right)\right)$ such that:

$$\left\| \Psi\left(t+j\delta,0\right)-\Psi^{*}\left(t+j\delta,0\right)\right\| \le \varepsilon_1\left(t+j\delta\right); \quad \left\| \tilde{A}\left(t+\delta_j\right)-\tilde{A}^{*}\right\| \le \varepsilon_2\left(t+j\delta\right) \tag{2.19}$$

$$\begin{aligned}\left\| \tilde{A}\left(t\right)\Psi\left(t,0\right)-\tilde{A}^{*}\Psi^{*}\left(t,0\right)\right\| &= \left\| \tilde{A}\left(t\right)\Psi\left(t,0\right)-\tilde{A}^{*}\Psi^{*}\left(t,0\right)\right\| \\ &\le \left\| \left(\tilde{A}\left(t\right)-A^{*}\right)\Psi\left(t,0\right)+A^{*}\left(\Psi\left(t,0\right)-\Psi^{*}\left(t,0\right)\right)\right\| \\ &\le \left\| A^{*}\left(\Psi\left(t,0\right)-\Psi^{*}\left(t,0\right)\right)\right\| + \left\| \left(\tilde{A}\left(t\right)-A^{*}\right)\Psi\left(t,0\right)\right\| \\ &\le \varepsilon_1\left(t+j\delta\right)\left\| A^{*}\right\| + \varepsilon_2\left(t+j\delta\right)\left\| \Psi\left(t+j\delta,0\right)\right\| \end{aligned} \tag{2.20}$$

$\forall t \ge max\left(t_1\left(\varepsilon_1\right), t_2\left(\varepsilon_2\right)\right)$. Thus, one gets from (2.20) into (2.19):

$$\left\| \dot{\Psi}\left(t+j\delta,0\right)-\dot{\Psi}^{*}\left(t+j\delta,0\right)\right\| \le \varepsilon_1\left(t+j\delta\right)\left(\left\| A\right\|+\left\| A^{*}\right\|\right)+\varepsilon_2\left(t+j\delta\right)\left\| \Psi\left(t+j\delta,0\right)\right\| \tag{2.21}$$

$\forall t \ge max\left(t_1\left(\varepsilon_1\right), t_2\left(\varepsilon_2\right)\right)$. One gets from (2.19) and (2.21) proceeding by complete induction:

$$\begin{aligned}&\left\| \Psi\left(t+(j+1)\delta,0\right)-\Psi^{*}\left(t+(j+1)\delta,0\right)\right\| \\ &\quad\le \varepsilon_1\left(t+j\delta\right)\left(\delta\left(\left\| A\right\|+\left\| A^{*}\right\|\right)+1\right)+\varepsilon_2\left(t+j\delta\right)\delta\left\| \Psi\left(t+j\delta,0\right)\right\| \\ &\quad\le \varepsilon_1\left(t+(j+1)\delta\right)\left(\left(\left\| A\right\|+\left\| A^{*}\right\|\right)+1\right)+\varepsilon_2\left(t+(j+1)\delta\right)\left\| \Psi\left(t+(j+1)\delta,0\right)\right\| \end{aligned} \tag{2.22}$$

provided that

$$0<\delta<min\left(\frac{\varepsilon_1\left(t+(j+1)\delta\right)\left(\left(\left\| A\right\|+\left\| A^{*}\right\|\right)+1\right)-\varepsilon_1\left(t+j\delta\right)}{\left\| A\right\|+\left\| A^{*}\right\|}, \frac{\varepsilon_2\left(t+(j+1)\delta\right)\left\| \Psi\left(t+(j+1)\delta,0\right)\right\|}{\varepsilon_2\left(t+j\delta\right)\left\| \Psi\left(t+j\delta,0\right)\right\|}\right)$$

Thus,

$$\lim_{N_0 \ni j \to \infty}\left\| \Psi\left(t+(j+1)\delta,0\right)-\Psi^{*}\left(t+(j+1)\delta,0\right)\right\| = \lim_{S_j\left(t,\delta,\varepsilon_1,\varepsilon_2\right) \ni t_j \to \infty}\left\| \Psi\left(t+(j+1)\delta,0\right)-\Psi^{*}\left(t+(j+1)\delta,0\right)\right\| = 0$$

where $S_j\left(t,\delta,\varepsilon_1,\varepsilon_2\right) := \left\{t_j := t+j\delta : j \in \boldsymbol{N}_0\right\}$; $\forall t \ge t_0 := max\left(t_1\left(\varepsilon_1\right), t_2\left(\varepsilon_2\right)\right)$. As a result, for any $\delta \in \boldsymbol{R}_+$.

$$\lim_{t \to \infty}\left\| \Psi\left(t,0\right)-\Psi^{*}\left(t,0\right)\right\| = \lim_{\bigcup_{t \in [t_0, t_0+\delta)} S_j\left(t,\delta,\varepsilon_1,\varepsilon_2\right) \ni t_j \to \infty}\left\| \Psi\left(t+(j+1)\delta,0\right)-\Psi^{*}\left(t+(j+1)\delta,0\right)\right\| = 0 \tag{2.23}$$

Property (i) has been proven. The proof of Property (ii) is quite similar to that of Theorem 2.1 (iii) so that it is omitted. □

A direct consequence of Theorem 2.1(iii) and Theorem 2.2 (ii) is that, since the discrete-time system has a limiting equation if all the coefficients of the differential equation (2.1) converge asymptotically to finite limits for any given sampling period T, then there is also (equivalently) a limiting difference equation which describes the limiting solution of (2.1) at consecutive samples for any given sampling period T. The discrete-time equation associated with (2.1) is:

$$x_{k+1} = e_1^T \Phi_k w_k = \varphi_{1k}^T w_k = \sum_{j=1}^{n} \varphi_{1jk} x_{k+1-j} \tag{2.24}$$

$\forall k(\geq n) \in \boldsymbol{N}$, with arbitrary initial conditions $D^i x(0) = x_{0i} \in \boldsymbol{C}$ $(\forall i \in \overline{n-1} \cup \{0\})$ and some given sampling period $T \in \boldsymbol{R}_+$, where:

$$\Phi_k := (\varphi_{ijk}) = (\varphi_{1k}, \varphi_{2k}, \dots, \varphi_{nk})^T \ ; \ \varphi_{ik} := (\varphi_{i1k}, \varphi_{i2k}, \dots, \varphi_{ink})^T \tag{2.25a}$$

$$w_k := (x_k, x_{k-1}, \dots, x_{k+1-n})^T \tag{2.25b}$$

and its limiting version for $t = kT$ as $\boldsymbol{N} \ni k \to \infty$ is:

$$x_{k+1}^* = e_1^T \Phi^* w_k^* = \varphi_1^{*T} w_k^* = \sum_{j=1}^{n} \varphi_{1j}^* x_{k+1-j}^* \tag{2.26}$$

$\forall k(\geq n) \in \boldsymbol{N}$, with arbitrary initial conditions $D^i x(0) = x_{0i} \in \boldsymbol{C}$ $(\forall i \in \overline{n-1} \cup \{0\})$, where:

$$\Phi^* := (\varphi_{ij}^*) = (\varphi_1^*, \varphi_2^*, \dots, \varphi_n^*)^T \ ; \ \varphi_i^* := (\varphi_{i1}^*, \varphi_{i2}^*, \dots, \varphi_{in}^*)^T \tag{2.27a}$$

$$w_k^* := (x_k^*, x_{k-1}^*, \dots, x_{k+1-n}^*)^T \tag{2.27b}$$

A unique difference equation (2.24)-(2.25) exists for each given sampling period $T \in \boldsymbol{R}_+$. Also, a unique corresponding limit equation (2.26)-(2.27) exists for the same sampling period provided that $\alpha_i : \boldsymbol{R}_{0+} \to \boldsymbol{C}$ ; $\forall i \in \bar{n} \cup \{0\}$.

## 3. Main results about stability, error bounds and asymptotic expansions

Eqn. 2.1 is compared with its limiting counterpart and the corresponding error through the subsequent result by using the corresponding vector discrete-time equations (2.24)-(2.25) and (2.26)-(2.27):

THEOREM 3.1. Assume that $\lim_{t \to \infty} A(t) = A^*$ and define the error sequence at sampling instants $\tilde{w}_k := w_k - w_k^*$ as well as the parametrical error of the discrete-time difference system with respect to its limiting equation $\tilde{\Phi}_k := \Phi_k - \Phi^*$ for any given sampling period $T \in \boldsymbol{R}_+$. The following properties hold:

(i) If $\max_{0\le i\le j}\left\| w_{k+i}\right\| \le W<\infty$ for some $j\in \boldsymbol{N}_0$, then $\max_{0\le i<\infty}\left\| w_{k+i}\right\| \le W<\infty$ and $\limsup_{k\to\infty}\left\| w_{k+i}\right\| \le W<\infty$ for some $W\in \boldsymbol{R}_+$ provided that

$0<K^* \le \dfrac{1-\tilde{\rho}}{\rho^{*N}\left(1-\tilde{\rho}+\tilde{K}\tilde{\rho}^k\left(1-\tilde{\rho}^N\right)\max\left(1,\rho^{*j-1}\right)\right)}$ where $K^*$, $\tilde{K}$ (being norm-dependent) and $\rho^*$ and $\tilde{\rho}\in \boldsymbol{R}_+\cap[0,1)$ are real constants such that $\left\|\Phi^{*i}\right\|\le K^*\rho^{*i}$, $\left\|\prod_{\ell=k}^{k+i-1}\tilde{\Phi}_\ell\right\|\le \tilde{K}\tilde{\rho}^i$ (i.e. the limiting differential system is exponentially stable and the parametrical error with respect to the limiting discrete-time difference system converges exponentially fast to zero) ; $\forall i,k\in \boldsymbol{N}_0$.

(ii) $$\left|\left\|\tilde{w}_{k+j}\right\|-\left\|\tilde{w}_k\right\|\right| \le \frac{\tilde{K}\tilde{\rho}^k\left(1-\tilde{\rho}^j\right)}{1-\tilde{\rho}}\max_{0\le i\le j}\left\| w^*_{k+i}\right\|<\infty$$

$\forall k\in \boldsymbol{N}_0, \forall j(\ge j_0)\in \boldsymbol{N}_0$ and a sufficiently large finite $j_0\in \boldsymbol{N}_0$ provided that $\rho^*\in[0,1]$, $\rho\in[0,1)$ and $\tilde{\rho}\in[0,1)$. Also $\limsup_{k\to\infty}\left(\left\|\tilde{w}_{k+j}\right\|-\left\|\tilde{w}_k\right\|\right)=0$, $\forall j\in \boldsymbol{N}_0$, and $\limsup_{k\to\infty\ j\to\infty}\left(\left\|\tilde{w}_{k+j}\right\|-\left\|\tilde{w}_k\right\|\right)=0$.

(iii) Assume that $\tilde{\rho}\in\left[0,\rho^*\right)\subset[0,1)$ and $\rho^*\in\left[0,\rho_0/K^*\right)\in[0,1)$ for some $\rho_0\in \boldsymbol{R}_+$ such that $[0,\rho_0)\subset[0,1)$. Then,

$$\max_{0\le k<\infty}\left\|\tilde{w}_k\right\|\le\left(1-K^*\left(\rho^*+\frac{\tilde{K}\tilde{\rho}^k}{1-\tilde{\rho}/\rho^*}\right)\right)^{-1}\frac{K^*\tilde{K}\tilde{\rho}^k}{1-\tilde{\rho}/\rho^*}\max_{0\le k<\infty}\left\| w^*_k\right\|<\infty$$

and $\lim_{k\to\infty}\tilde{w}_k=\lim_{k\to\infty}w_k=\lim_{k\to\infty}w^*_k=0$.

*Proof*: (i) Note that

$$w_{k+1}=\Phi_k w_k=\Phi^* w_k+\tilde{\Phi}_k w_k \tag{3.1}$$

$$\begin{aligned}\left\| w_{k+j}\right\| &=\left\|\Phi^{*j}w_k+\sum_{i=0}^{j-1}\Phi^{*j-i-1}\tilde{\Phi}_{k+i}w_{k+i}\right\| \\ &\le K^*\rho^{*j}\left\| w_k\right\|+\sum_{i=0}^{j-1}K^*\rho^{*j-i-1}\tilde{K}\tilde{\rho}^{k+i}\left\| w_{k+i}\right\| \\ &\le K^*\left(\rho^{*j}\left\| w_k\right\|+\frac{\tilde{K}\tilde{\rho}^k\left(1-\tilde{\rho}^j\right)\max\left(1,\rho^{*j-1}\right)}{1-\tilde{\rho}}\max_{0\le i\le j}\left\| w_{k+i}\right\|\right)\end{aligned}$$

$$\leq K^* \rho^{*\,j} \left( 1 + \frac{\tilde{K}\,\tilde{\rho}^{\,k}\left(1-\tilde{\rho}^{\,j}\right) max\left(1,\ \rho^{*\,j-1}\right)}{1-\tilde{\rho}} \right) \max_{0\leq i\leq j} \left\| w_{k+i} \right\| \tag{3.2}$$

Note from (3.1) that

$$0 < K^* \leq \frac{1-\tilde{\rho}}{\rho^{*\,j}\left(1-\tilde{\rho}+\tilde{K}\,\tilde{\rho}^{\,k}\left(1-\tilde{\rho}^{\,j}\right) max\left(1,\ \rho^{*\,j-1}\right)\right)j} \Rightarrow K^*\left( \rho^{*\,j} + \frac{\tilde{K}\,\tilde{\rho}^{\,k}\left(1-\tilde{\rho}^{\,j}\right) max\left(1,\ \rho^{*\,j-1}\right)}{1-\tilde{\rho}} \right) \leq 1$$

$$\Rightarrow \max_{0\leq i\leq j+1} \left\| w_{k+i} \right\| \leq \max_{0\leq i\leq j} \left\| w_{k+i} \right\|$$

and it implies also that $0 < K^* \leq \frac{1-\tilde{\rho}}{\rho^{*\,j}\left(1-\tilde{\rho}+\tilde{K}\,\tilde{\rho}^{\,k+1}\left(1-\tilde{\rho}^{\,j}\right) max\left(1,\ \rho^{*\,j-1}\right)\right)}$

Since $0 < \tilde{\rho}^{\,k} \leq \tilde{\rho}^{\,k+1} < 1$ and the proof that $\max_{0\leq i<\infty} \left\| w_{k+i} \right\| \leq W < \infty$ and $\limsup_{k\to\infty} \left\| w_{k+i} \right\| \leq W < \infty$ then follow by complete induction. Then, property (i) has been proven.

(ii) Note that (3.1) is equivalent to the identities below:

$$\tilde{w}_{k+1} = \Phi_k w_k - \Phi^* w_k^* = \Phi_k \tilde{w}_k + \tilde{\Phi}_k w_k^* = \Phi^* \tilde{w}_k + \tilde{\Phi}_k w_k \tag{3.3}$$

so that, since $\left\| \Phi^{\,i} \right\| \leq K \rho^{\,i}$ for some real constants $K \in \mathbf{R}_+$ and $\rho \in (0,1)$; $\forall i \in \mathbf{N}_0$,

$$\left\| \tilde{w}_{k+j} \right\| = \left\| \Phi^{\,j}\, \tilde{w}_k + \sum_{i=0}^{j-1} \Phi^{\,j-i-1} \tilde{\Phi}_{k+i}\, w_{k+i}^* \right\|$$

$$\leq K\rho^{\,j} \left\| \tilde{w}_k \right\| + \sum_{i=0}^{j-1} K\rho^{\,j-i-1}\, \tilde{K}\,\tilde{\rho}^{\,k+i} \left\| w_{k+i}^* \right\|$$

$$\leq K\left( \rho^{\,j} \left\| \tilde{w}_k \right\| + \frac{\tilde{K}\,\tilde{\rho}^{\,k}\left(1-\tilde{\rho}^{\,j}\right)}{1-\tilde{\rho}} \max_{0\leq i\leq j} \left\| w_{k+i}^* \right\| \right)$$

$$\Rightarrow \left| \left\| \tilde{w}_{k+j} \right\| - \left\| \tilde{w}_k \right\| \right| \leq \left\| \tilde{w}_{k+j} \right\| - K\rho^{\,j} \left\| \tilde{w}_k \right\| \leq \frac{\tilde{K}\,\tilde{\rho}^{\,k}\left(1-\tilde{\rho}^{\,j}\right)}{1-\tilde{\rho}} \max_{0\leq i\leq j} \left\| w_{k+i}^* \right\| < \infty \tag{3.4}$$

$\forall k \in \mathbf{N}_0, \forall j\left(\geq j_0\right) \in \mathbf{N}_0$ and a sufficiently large finite $j_0 \in \mathbf{N}_0$ since $\left\| w_k^* \right\| < \infty$; $\forall k \in \mathbf{N}_0$ since $\rho^* \in [0,1]$, $\rho \in [0,1)$ and $\tilde{\rho} \in [0,1)$; $\forall k \in \mathbf{N}_0$. Furthermore, $\limsup_{k\to\infty} \left( \left\| \tilde{w}_{k+j} \right\| - \left\| \tilde{w}_k \right\| \right) = 0$, $\forall j \in \mathbf{N}_0$, and $\limsup_{k\to\infty\ j\to\infty} \left( \left\| \tilde{w}_{k+j} \right\| - \left\| \tilde{w}_k \right\| \right) = 0$. Property (ii) has been proven.

(iii) From the last expression in (3.2):

$$\left\| \tilde{w}_{k+j} \right\| = \left\| \Phi^{*\,j}\, \tilde{w}_k + \sum_{i=0}^{j-1} \Phi^{*\,j-i-1}\, \tilde{\Phi}_{k+i}\, w_{k+i} \right\|$$

$$\le K^* \left( \rho^{*\,j} \left\| \tilde{w}_k \right\| + \tilde{K} \rho^{*\,j-1} \tilde{\rho}^{\,k} \sum_{i=0}^{j-1} \left( \frac{\tilde{\rho}}{\rho^*} \right)^i \left\| w_{k+i} \right\| \right)$$

$$\le K^* \left( \rho^{*\,j} \left\| \tilde{w}_k \right\| + \tilde{K} \rho^{*\,j-1} \tilde{\rho}^{\,k} \frac{1-\left( \tilde{\rho} / \rho^* \right)^j}{1-\tilde{\rho}/\rho^*} \max_{0 \le i \le j} \left\| w_{k+i} \right\| \right)$$

$$\le K^* \rho^{*\,j} \left( \left\| \tilde{w}_k \right\| + \tilde{K} \rho^{*\,-1} \tilde{\rho}^{\,k} \frac{1-\left( \tilde{\rho} / \rho^* \right)^j}{1-\tilde{\rho}/\rho^*} \left( \max_{0 \le i \le j} \left\| w^*_{k+i} \right\| + \max_{0 \le i \le j} \left\| \tilde{w}_{k+i} \right\| \right) \right)$$

$$\le K^* \rho^{*\,j-1} \left( \rho^* + \tilde{K} \tilde{\rho}^{\,k} \frac{1-\left( \tilde{\rho} / \rho^* \right)^j}{1-\tilde{\rho}/\rho^*} \right) \max_{0 \le i \le j} \left\| \tilde{w}_{k+i} \right\|$$

$$+ K^* \tilde{K} \rho^{*\,j} \tilde{\rho}^{\,k} \frac{1-\left( \tilde{\rho} / \rho^* \right)^j}{1-\tilde{\rho}/\rho^*} \max_{0 \le i \le j} \left\| w^*_{k+i} \right\|$$

$$\le K^* \left( \rho^* + \frac{\tilde{K} \tilde{\rho}^{\,k}}{1-\tilde{\rho}/\rho^*} \right) \max_{0 \le i \le j} \left\| \tilde{w}_{k+i} \right\| + \frac{K^* \tilde{K} \tilde{\rho}^{\,k}}{1-\tilde{\rho}/\rho^*} \max_{0 \le i \le j} \left\| w^*_{k+i} \right\| ; \quad \forall j \in \mathbf{N} \quad (3.5)$$

Comparing the first and last terms of (3.5), it is obvious that

$$\max_{0 \le i \le j} \left\| \tilde{w}_{k+i} \right\| \le K^* \left( \rho^* + \frac{\tilde{K} \tilde{\rho}^{\,k}}{1-\tilde{\rho}/\rho^*} \right) \max_{0 \le i \le j} \left\| \tilde{w}_{k+i} \right\| + \frac{K^* \tilde{K} \tilde{\rho}^{\,k}}{1-\tilde{\rho}/\rho^*} \max_{0 \le i \le j} \left\| w^*_{k+i} \right\|$$

$$; \quad \forall j \in \mathbf{N} \qquad (3.6)$$

Since $K^* \rho^* \in [0, \rho_0) \subset [0,1)$ and $\tilde{\rho} \in \lfloor 0, \rho^* ) \subset [0,1)$, $1 > K^* \left( \rho^* + \tilde{K} \tilde{\rho}^{\,k} \right)$ for some sufficiently large and finite $k_0 \in \mathbf{N}_0$ and $\forall k (\ge k_0) \in \mathbf{N}_0$. Also, $\max_{0 \le k < \infty} \left\| w^*_k \right\| < \infty$ since $\rho^* \in \lfloor 0, \rho_0 / K^* ) \in [0,1)$ so that Property (iii) follows from (3.6) and $\lim_{k \to \infty} \max_{0 \le i < \infty} \left\| \tilde{w}_{k+i} \right\| \le \left( \frac{1}{1-K^* \rho^*} \right) \frac{K^* \tilde{K}}{1-\tilde{\rho}/\rho^*} \max_{0 \le k < \infty} \left\| w^*_{k+i} \right\| \left( \lim_{k \to \infty} \tilde{\rho}^{\,k} \right) = 0$. □

Again

$$\tilde{w}_{k+1} = \Phi_k w_k - \Phi^* w^*_k = \Phi_k \tilde{w}_k + \tilde{\Phi}_k w^*_k = \Phi^* \tilde{w}_k + \tilde{\Phi}_k w_k \qquad (3.7)$$

so that one gets from the left-hand –side and the third identity of (3.7) :

$$\limsup_{k\to\infty}\left|\frac{\left\|\tilde{w}_{k+j}\right\|-\left\|\Phi^{*^{j}}\tilde{w}_k\right\|}{\max_{0\le i\le j}\left\|w_{k+i}\right\|}\right|\le\frac{\left(1-\left(\tilde{\rho}/\rho^*\right)^j\right)K^*\tilde{K}}{1-\tilde{\rho}/\rho^*}\limsup_{k\to\infty}\left(\rho^{*^{j-1}}\tilde{\rho}^k\right)=0$$

provided that $0\le\tilde{\rho}<\min\left(1,\rho^*\right)$ for $\rho^*\in\boldsymbol{R}_+$, $\forall j\in\boldsymbol{N}$. Also, since

$$\left\|\tilde{w}_{k+j}\right\|\le K^*\left(\rho^{*^{j}}\left\|\tilde{w}_k\right\|+\tilde{K}\rho^{*^{j-1}}\tilde{\rho}^k\sum_{i=0}^{j-1}\left(\frac{\tilde{\rho}}{\rho^*}\right)^i\left\|w_{k+i}\right\|\right)$$

$$\le K^*\rho^{*^{j}}\left(\left\|\tilde{w}_k\right\|+\frac{\tilde{K}\tilde{\rho}^k}{\rho^*}\sum_{i=0}^{j-1}\left(\frac{\tilde{\rho}}{\rho^*}\right)^i\left(\left\|\tilde{w}_{k+i}\right\|+\left\|w^*_{k+i}\right\|\right)\right)$$

$$\limsup_{k\to\infty}\left|\frac{\left\|\tilde{w}_{k+j}\right\|-\left\|\Phi^{*^{j}}\tilde{w}_k\right\|}{\max_{0\le i\le j}\left\|w_{k+i}\right\|}\right|\le\frac{\left(1-\left(\tilde{\rho}/\rho^*\right)^j\right)K^*\tilde{K}}{1-\tilde{\rho}/\rho^*}\limsup_{k\to\infty}\left(\rho^{*^{j-1}}\tilde{\rho}^k\right)=0 \quad \square$$

The following result follows directly from Theorem 3.3 and (2.25)-(2.28) concerning with the relationships of the solutions the differential equation (2.1) with that of its limiting equation defined for $\alpha_i^*=\lim_{t\to\infty}\alpha_i(t)$, $\forall i\in\overline{n}\cup\{0\}$ obtained from the associated differential system (2.2)-(2.3) and its limiting version since the results in Theorem 3.1 hold for any given sampling period T.

COROLLARY 3.2. Consider the time-varying differential equation (2.1) together with its limiting counterpart $\sum_{i=0}^{n}\alpha^*_{n-i}D^ix^*(t)=0$, $\quad D^ix(0)=D^ix^*(0)=x_{0i}\in\boldsymbol{C}\quad(\forall i\in\overline{n-1}\cup\{0\})$ defined for $A^*=\lim_{t\to\infty}A(t)$; i.e. $\alpha_i^*=\lim_{t\to\infty}\alpha_i(t)$, $\forall i\in\overline{n-1}\cup\{0\}$. Thus, the following properties hold:

(i) The differential equation (2.1) and the differential system (2.2)-(2.3) satisfy the properties
$\left|\left\|\tilde{z}(t+\tau)\right\|-\left\|\tilde{z}(t)\right\|\right|<\infty$ ; $\left|\tilde{x}(t+\tau)-\tilde{x}(t)\right|<\infty$
$\forall\tau\in\boldsymbol{R}_{0+}$ and $\forall t\in\boldsymbol{R}_{0+}$ being finite and provided that the limiting differential equation $\sum_{i=0}^{n}\alpha^*_{n-i}D^ix^*(t)=0$ is exponentially stable with stability abscissa $-\rho^*_\alpha<0$ and $\alpha_i(t)=\alpha_i^*\left(1-e^{-\alpha(t-t_0)}\right)$; $\forall t\ge t_0$ ; $\forall i\in\overline{n-1}\cup\{0\}$ (for some finite $t_0\in\boldsymbol{R}_{0+}$) and some $\alpha\in\boldsymbol{R}_+$.

(ii) Assume that $-\alpha < -\rho_{\alpha}^{*} < -\rho_0 / K^{*} < 0$ for some $\rho_0 \in \boldsymbol{R}_+$. Then, $\lim_{t\to\infty} \tilde{x}(t) = \lim_{t\to\infty} x(t) = \lim_{t\to\infty} x^{*}(t) = 0$ and $\lim_{t\to\infty} \tilde{z}(t) = \lim_{t\to\infty} z(t) = \lim_{t\to\infty} z^{*}(t) = 0$.

*Proof*: It is a direct consequence of Theorem 3.1 (ii)-(iii) since the existence of real numbers $\rho_T^{*} \in [0,1]$, $\rho_T \in [0,1)$ and $\tilde{\rho}_T \in [0,1)$, depending on $T \in \boldsymbol{R}_+$, leading to the Properties (ii)-(iii) of Theorem 3.1 for any given sampling period $T \in \boldsymbol{R}_+$, implies the existence $\rho^{*} := \min_{T\in(0,T^{*}]} \rho_T^{*} \in [0,1]$, $\rho := \min_{T\in(0,T^{*}]} \rho_T \in [0,1]$ and $\tilde{\rho} := \min_{T\in(0,T^{*}]} \tilde{\rho}_T \in [0,1]$ which lead to the properties for all $T \in (0,T^{*}]$ for any given arbitrary $T^{*} \in \boldsymbol{R}_+$. Thus, the Properties (i)-(ii) of the continuous-time differential equation and differential system follow. □

The following result relates the boundedness and unboundedness of the solutions of (2.1) with their limiting equation counterparts depending on the modulus of a dominant characteristic root of the limiting equation by extending a parallel previous result for the discrete-time case, [1] . In the sequel, rename the solution sequence $\{x_k\}_0^{\infty}$ satisfying $x_{k+n} = e_1^T z_{k+n}$, obtained from the discrete system (2.11), as $\{x_{Tk}\}_0^{\infty}$ to reflect its dependence on the sampling period T for purposes of comparison of the solution sequences for distinct values of T.

THEOREM 3.3. Assume that $\lim_{t\to\infty} A(t) = A^{*}$ and that $x : \boldsymbol{R}_{0+} \to \boldsymbol{C}$ is a nontrivial solution of the (2.1) any set of given nonzero initial conditions. Then, the following properties hold:

(i) $|x(kT)| \le \limsup_{\boldsymbol{N}_0 \ni k\to\infty} |x(kT)| + \varepsilon \le \mu_T^k + \varepsilon$ for any sampling period $T \in \boldsymbol{R}_+$, $\forall k \ge k_0$, for any given arbitrarily small $\varepsilon \in \boldsymbol{R}_+$, and some finite $k_0 \in \boldsymbol{N}_0$ depending on $\varepsilon$, where: $\Lambda_T = \Lambda(\mu_T) := \{ z \in \boldsymbol{C} : \Delta_T(z) = 0, |z| = \mu_T \} \ne \varnothing$ ,

$\Delta_T(z) := \sum_{i=1}^{n} \psi_{iT}^{*} z^{k-i}$ is the characteristic polynomial of $\Psi^{*}(T,0) = \left(\psi_{ij}^{*}(T,0)\right) := e^{nA^{*}T}$, and

$a_{iT}^{*}$ are the entries of its last row with z being the Z- transform argument .

(ii) Assume that the function of the sampling period $\mu_T := \limsup_{k\to\infty} \sqrt[k]{|x(kT)|} : \boldsymbol{R}_{0+} \to \boldsymbol{R}_{0+}$ is such that $\mu_T \le 1$ for some given $T \in \boldsymbol{R}_{0+}$. Then , the difference Poincaré system (2.11) and its limiting counterpart are both globally stable in the Lyapunov´s sense. Also, no characteristic root of the limiting difference equation associated with (2.11)-(2.12) has modulus greater than unity for any arbitrary sampling period $T \in \boldsymbol{R}_{0+}$.

(iii) If $\mu_T < 1$ for some $T \in \boldsymbol{R}_+$ then both the Poincaré system (2.11)-(2.12) and its associated limiting difference system are both globally asymptotically stable for the sampling period $T \in \boldsymbol{R}_+$.

(iv) If $\mu_T > 1$ for some $T \in \boldsymbol{R}_+$ then the difference Poincaré system (2.11)-(2.12), its limiting counterpart and its respective associate limiting equations are unstable for any sampling period and so it is the continuous-time differential equation (2.1).

*Proof*: For any sampling period $T \in \boldsymbol{R}_+$, the discrete system (2.11)-(2.12) satisfies $\lim_{k \to \infty} \Psi((k+n)T, kT) = e^{A^* nT}$ and $x_k = x(kT) = e_1^T z(kT) = e_1^T z_k$, the solution of (2.11)-(2.12), is also the solution of (2.1) at sampling instants for corresponding related initial conditions. The solution of the limiting equation of the discrete system (2.11)-(2.12) is $z^*((k+1)T) := z_{k+1}^* = e^{nA_T^* T} z_k^*$ and, equivalently, that of the discrete limiting equation associated with (2.1) is $x_k^* := x^*(kT) = e_1^T z_k^*$, obtained from (2.11)-(2.12) with the replacement $\Psi((k+n)T, kT) \to \Psi^*((k+n)T, kT) = e^{A^* nT}$, equivalent to the replacement $\Psi((k+1)T, kT) \to \Psi^*(T, 0) = e^{A^* T}$ which is a constant matrix independent of k for each given T. From Theorem 2.1 in [1], $|x(kT)| \le \limsup_{\boldsymbol{N}_0 \ni k \to \infty} |x(kT)| + \varepsilon \le \mu_T^k + \varepsilon$ for any given arbitrarily small $\varepsilon \in \boldsymbol{R}_+$, $\forall k \ge k_0$ and some finite $k_0 \in \boldsymbol{N}_0$ depending on $\varepsilon$, where $\mu_T := \limsup_{k \to \infty} \sqrt[k]{|x(kT)|}$ so that $\Lambda_T \ne \emptyset, \forall T \in \boldsymbol{R}_+$. The above property implies and it is implied by $\mu_T$ being one of the characteristic roots of the limiting equation $x_{Tk}^* = x^*(kT) = e_1^T z_k^*$. Now, it follows directly that:

a) If $\mu_T \le 1$ then the difference Poincaré system (2.11) and its limiting counterpart are both globally stable in the Lyapunov´s sense since $|x(kT)| \le 1 + \varepsilon < \infty$ for any given arbitrarily small $\varepsilon \in \boldsymbol{R}_+$, $\forall k \ge k_0$ and some finite $k_0 \in \boldsymbol{N}_0$ depending on $\varepsilon$ and $|x(kT)| \le M < \infty$ for some finite $M \in \boldsymbol{R}_+$, depending on the initial conditions, $\forall k (\in \boldsymbol{N}_0) < k_0$. Since the solution of the solution of the differential Poincaré system for corresponding initial conditions coincides at sampling instants with that of the difference system , i.e. $x(kT) = x_k$ and the solution of (2.1) is a continuously- differentiable real function on $\boldsymbol{R}_{0+}$, it follows that it cannot be unbounded in-between consecutive sampled values where it is bounded .Thus, $|x(t)| \le K_1 < \infty$, $\forall t \in [kT, (k+1)T)$, $\forall k \in \boldsymbol{N}_0$. Its limiting counterpart is also uniformly bounded on $\boldsymbol{R}_{0+}$. Also, $\mu_{T_s} \le 1$ for any other $T_s (\ne T) \in \boldsymbol{R}_+$ since $x(kT_s)$ is the solution of (2.1) at $t = kT_s$ for identical initial conditions and it is continuously differentiable and bounded at sampling instants so that it is uniformly bounded on $\boldsymbol{R}_{0+}$. Note finally, that all the roots of the characteristic limiting equation associated with the Poincaré differential system (2.11)-(2.12) have modulus less than unity since such a system is Lyapunov´s stable. Property (i) has been proven.

b) It turns out from the above stability property that if $\mu_T \le 1$ for some sampling period T then no characteristic root of the limiting difference equation associated with the limiting system (2.11)-(2.12) has modulus greater than unity for any sampling period. Property 8ii) has been proven.

c) If $\mu_T < 1$ for some $T \in \boldsymbol{R}_+$ then both the Poincaré system (2.11)-(2.12) and its associated limiting difference system are both globally asymptotically stable for the sampling period $T \in \boldsymbol{R}_+$ since $\left| x(kT) \right|$ is uniformly bounded, $\forall\, k \in \boldsymbol{N}_0$ and $\lim_{\boldsymbol{N}_0 \ni k \to \infty} x(kT) = \lim_{\boldsymbol{N}_0 \ni k \to \infty} \mu_T^k = 0$. Also, the limiting differential equation associated with (2.1), the limiting difference system associated with (2.11)-(2.12) and the associated limiting difference equation are globally asymptotically stable for any other sampling period since $\lim_{k \to \infty} x_{Tk}^* = \lim_{k \to \infty} \mu_T^k = 0 \Rightarrow \lim_{t \to \infty} x^*(t) = 0$ ( since the limiting difference equation associated with (2.11)-(2.12) and the limiting differential system associated with (2.1) are linear time-invariant . Then, $\lim_{\boldsymbol{N}_0 \ni k \to \infty} x^*\left(kT^{\prime}\right) = 0$, $\forall T^{\prime} \in \boldsymbol{R}_+$. Property (iii) has been proven,

d) If $\mu_T > 1$ for some $T \in \boldsymbol{R}_+$ then the difference Poincaré system (2.11)-(2.12), its limiting counterpart and its respective associate limiting equations are unstable. The property also holds for any $T_s(\ne T) \in \boldsymbol{R}_+$ and for the continuous-time differential equation (2.1). Property (iv) has been proven. □

REMARK 3.4. Note from Theorem 3.3 with $\mu_T < 1$ for sampling period T implies the following issues. The differential limiting equation associated with (2.1) and the differential Poincaré limiting system, obtained from (2.2) by replacing $A(t)$ with $A^*$, are globally asymptotically stable. As a result, all the Poincaré limiting difference equations and their corresponding limiting difference systems obtained from (2.1) and (2.11)-(2.12), equivalently the limiting difference equation (2.26)-(2.27) and its corresponding limiting difference system, with the replacement $A(t) \to A^*$ are globally asymptotically stable for any sampling period. However, the test $\mu_T < 1$ does not guarantee the global asymptotic stability of the Poincaré differential equation (2.1), since it is not time-invariant, but only that $\lim_{\boldsymbol{N}_0 \ni k \to \infty} x(kT) = 0$ and $x(t)$ is uniformly bounded on $\boldsymbol{R}_{0+}$. □

Now, the characteristic solutions of the limiting differential equation of (2.1) obtained from (2.1) by replacing $\alpha_i$ with $\alpha_i^* = \lim_{t \to \infty} \alpha_i(t)$, $\forall i \in \overline{n-1} \cup \{0\}$ are obtained for all time by applying calculus of residues to the discrete limiting equation (2.16)-(2.27) by defining time-varying sampling periods associated to each time instant where the computation is performed.

THEOREM 3.5. Let $\lambda \in \boldsymbol{C}$ be some eigenvalue of $A^{*} = \left( a_{ij}^{*} \right)$ of multiplicity $m_{\lambda}$. For any given (finite or infinite) $t \in \boldsymbol{R}_{0+}$, define $T_t := \frac{t}{k_t} \in \boldsymbol{R}_{0+}$ for some $k_t \in \boldsymbol{N}_0$, being both possibly dependent on t, such that either $T_t < \infty$ for any given $\lambda \in \boldsymbol{C}$ or $T_t \le \infty$, for $\lambda \in \boldsymbol{C}_{0+}$. Let $\Delta_T(z) := \sum_{i=1}^{n} \psi_{iT_t}^{*} z^{n-i}$ be the characteristic polynomial of $\Psi^{*}(T_t, 0) = e^{A^{*} T_t}$. Thus, the characteristic solution corresponding to $\lambda$ of the limiting differential equation of (2.1) may be calculated as:

$$x^{*}(t) = Res\left( z^{k_t - 1} \Delta_{T_t}^{-1}(z) f(z); \lambda \right) = Res\left( \frac{z^{k_t - 1} f(z)}{(z-\lambda)^{m_\lambda} q(z)}; \lambda \right) = c_{-1}\left( \frac{z^{k_t - 1} f(z)}{(z-\lambda)^{m_\lambda} q(z)}; \lambda \right) = p(k_t) \lambda^{k_t} \tag{3.8}$$

$\forall\, t \in \boldsymbol{R}_{0+}$ with $Res\, g_t(z;\lambda)$ denoting the residue of the complex valued function $g_t$ at a pole $\lambda$ of multiplicity $m_\lambda$ of $g_t$, $f : \boldsymbol{C} \to \boldsymbol{C}$ being any holomorphic function in a neighborhod of $z = \lambda$, $q(z)$ being any polynomial fulfilling $q(\lambda) \ne 0$, $p(k_t)$ being a polynomial function from $\boldsymbol{R}_{0+}$ to $\boldsymbol{R}_{0+}$ of degree less than $m_\lambda$, and $c_{-1}(g(z;\lambda))$ being the coefficient of the term $(z-\lambda)^{-1}$ of the Laurent series $\sum_{j=-m_\lambda}^{\infty} c_j(g_t(z;\lambda))(z-\lambda)^{-j}$ of the complex-valued function $g(z)$ below which is meromorphic in a region of the complex plane and which is identified with its Laurent series:

$$g_t(z) := z^{k_t - 1} \Delta_{T_t}^{-1}(z) f(z) = \frac{z^{k_t - 1} f(z)}{(z-\lambda)^{m_\lambda} q(z)} = \sum_{j=-m_\lambda}^{\infty} c_j(g_t(z;\lambda))(z-\lambda)^{-j} \tag{3.9}$$

*Proof*: If $\boldsymbol{R}_{0+} \ni T_t < \infty$ then $\lambda_{\Psi^*} := e^{\lambda T_t} \ne 0$ is a nonzero characteristic root of $e^{A^{*} T_t}$, $\forall \lambda \in \boldsymbol{C}$. If $\boldsymbol{R}_{0+} \ni T_t \le \infty$ then $\lambda_{\Psi^*} := e^{\lambda T_t} \ne 0$ if. In both cases, $e^{\lambda T_t}$ is always nonzero characteristic root of $e^{A^{*} T_t}$. Since $m_\lambda$ is the multiplicity of the eigenvalue $\lambda$ of $A^{*}$ then it is also the multiplicity of the nonzero eigenvalue $\lambda_{\Psi^*} = e^{\lambda T_t}$ of the fundamental matrix $e^{A^{*} T_t}$ of the solution of the associate limiting difference system $z_{k+1,t}^{*} = e^{A^{*} T_t} z_{kt}^{*}$. It follows by direct calculus of residues that (3.8) is true which is identical through (3.9) to

$$x^{*}(t) = \frac{1}{(m_\lambda - 1)!} \left. \frac{d^{m_\lambda - 1}}{dz^{m_\lambda - 1}} \right|_{z=\lambda} \left( (z-\lambda)^{m_\lambda} g_t(z) \right) = \frac{1}{(m_\lambda - 1)!} \left. \frac{d^{m_\lambda - 1}}{dz^{m_\lambda - 1}} \right|_{z=\lambda} \left( z^{k_t - 1} h(z) \right) \tag{3.10}$$

$\forall\, t \in \boldsymbol{R}_{0+}$, where $h(z) := p(z)/q(z)$. By direct computation of $\left. \frac{d^{m_\lambda - 1}}{dz^{m_\lambda - 1}} \right|_{z=\lambda} \left( z^{k_t - 1} h(z) \right)$, it follows that :

$$max\ degree\left(p\left(k_t\right)\right) < k_t - \left(k_t - 1 - pole\ /\ zero\ excess\left(h^{\left(m_\lambda - 1\right)}\right)\right) \le m_\lambda$$

(See [1] for the details about the above computation for Poincaré difference equations). □

REMARK 3.6. Note that Theorem 3.5 is valid to calculate at any time the characteristic solutions of a continuous-time time-varying linear differential equation with piecewise continuous-time coefficients with no zero characteristic root which possesses a limiting. The computation is performed through residues computation theory on an associate difference equation counterpart. In this way, the solution of the differential equation may be calculated by linear superposition of all the characteristic solutions without needing the explicit calculation of the solution via theory of solutions of differential equations. Note also from the statement of Theorem 3.5 that in order that such a result holds:

1) Either $T_t < \infty$, what implies $\boldsymbol{N}_0 \ni k_t \to \infty$ if $\boldsymbol{R}_{0+} \ni t \to \infty$ for arbitrary $\lambda \in \boldsymbol{C}$ being some eigenvalue of $A^*$, or
2) $T_t \le \infty$, for arbitrary $\lambda \in \boldsymbol{C}$ being some eigenvalue of $A^*$ with $Re\,\lambda \ge 0$, what implies that $k_t$ may be chosen so that either $\boldsymbol{N}_0 \ni k_t \to \infty$ if $\boldsymbol{R}_{0+} \ni t \to \infty$, as above, or fixed to any finite value $\boldsymbol{N}_0 \ni k_t < \infty$ if $\boldsymbol{R}_{0+} \ni t \to \infty$ implying that $\boldsymbol{R}_{0+} \ni T_t = \frac{t}{k_t} \to \infty$. □

The main result is the subsequent theorem which uses previous technical results from Theorems 3.3 and 3.5. The following result uses the standard "small -o" and "big- O" Landau´s notations:

THEOREM 3.7. Assume that the coefficient functions of (2.1) converge exponentially fast to their corresponding constant coefficients of the associate limiting equation; i.e. $\alpha_j(t) = \alpha_j^* + o\left(e^{-\sigma t}\right)$ ; $\forall\, j \in \bar{n} \cup \{0\}$ as $t \to \infty$ for some $\sigma \in \boldsymbol{R}_+$. If $x : \boldsymbol{R}_{0+} \to \boldsymbol{R}$ is a solution of (2.1), then either $x(t) = 0$ for all large $t \in \boldsymbol{R}_{0+}$, or one has the asymptotic expansion

$x(k_t T_t) = x^*(k_t T_t) + O\left(\left(\mu_{T_t} - \varepsilon_{T_t}\right)^{k_t}\right)$ as $\boldsymbol{R}_0 \ni t = k_t T_t \to \infty$ for any prefixed finite $T_t \in \boldsymbol{R}_+$, $\boldsymbol{N}_0 \ni k_t \to \infty$, and for some $\varepsilon_{T_t} \in \left(0, \mu_{T_t}\right)$ and some $\mu_{T_t} \in \boldsymbol{R}_+$ such that :

$$\Lambda_{T_t} = \Lambda\left(\mu_{T_t}\right) := \left\{ z \in \boldsymbol{C} : \Delta_{T_t}(z) = 0, |z| = \mu_{T_t} \right\} \ne \varnothing$$

where:

- $x^*(k_t T_t)$ is a nontrivial characteristic solution of the limiting equation (2.26)-(2.27) at large time $t = k_t T_t$ for sampling period $T_t$ which is identical to a nontrivial characteristic solution $x^*(t)$ of the limiting equation of (2.1), i.e. $\sum_{i=0}^{n} \alpha_{n-i}^* D^i x(t) = 0$ corresponding to the set $\Lambda_{T_t} = \Lambda\left(\mu_{T_t}\right)$

- $\Delta_{T_t}(z) := \sum_{i=1}^{n} \psi_{iT_t}^* z^{n-i}$ is the characteristic polynomial of $\Psi^*(T_t, 0) = e^{A^* T_t}$ .

*Proof*: First note that the coefficient $\varphi_{1n}^*$ in the limiting discrete- time equation (2.26)-(2.27) is nonzero irrespective of the sampling period. Otherwise, such a limiting equation would have a zero characteristic root. This is impossible since the eigenvalues of the fundamental matrix $e^{A^* t}$ of the limiting differential system are all nonzero for all $t \in \boldsymbol{R}_{0+}$. As a result $\mu_{T_t} > 0$ for any arbitrary $T_t \in \boldsymbol{R}_+$. Note from Theorem 3.3 that

$$\Lambda_{T_t} = \Lambda\left(\mu_{T_t}\right) := \left\{ z \in \boldsymbol{C} : \Delta_{T_t}(z) = 0, |z| = \mu_{T_t} \right\} \neq \varnothing$$

$$\Rightarrow \Lambda_\gamma = \Lambda(\gamma) := \left\{ s \in \boldsymbol{C} : \Delta(\gamma) = 0,\ \mu_{T_t} = p_\gamma\left(e^{-\gamma T_t}\right) \right\} \neq \varnothing \tag{3.11}$$

for some nonzero polynomial $p_\gamma\left(e^{-\gamma T_t}\right)$. Note that the differential equation (2.1) may be rewritten using its limiting equation as:

$$\sum_{i=0}^{n} \alpha_{n-i}^* D^i x(t) = s(t) \quad ; \quad s(t) := \sum_{i=0}^{n} \left(\alpha_{n-i}^* - \alpha_{n-i}(t)\right) D^i x(t) \tag{3.12}$$

Rewrite the difference equation (2.24) obtained from (2.1) for $t = k_t T_t$ and sampling period $T_t$ as follows:

$$x(t) \equiv x(k_t T_t) = \sum_{j=1}^{n} \varphi_{1j,k_t-1} x_{k_t-j} \equiv \sum_{j=1}^{n} \varphi_{1j,k_t-1}\left((k_t-1)T_t\right) x\left((k_t-j)T_t\right) \tag{3.13}$$

with the n- real vector $\varphi_{1j,k_t-1}$ being defined in (25.a). Note that x (t) is the solution of (2.1) for the given initial conditions $D^i x(0) = x_{0i} \in \boldsymbol{C}$ , $\forall i \in \overline{n-1} \cup \{0\}$ for $t = k_t T_t$, $\forall k_t \geq n$ provided that the initial conditions of (3.13) at $t_i = iT_t$ are $x_i \equiv x(iT_t)$ being the solution of (2.1) , $\forall i \in \overline{n-1} \cup \{0\}$. In the same way, the limiting difference equation for sampling period $T_t$ is:

$$x^*(t) \equiv x^*(k_t T_t) = \sum_{j=1}^{n} \varphi_{1j,k_t-1}^* x_{k_t-j} \equiv \sum_{j=1}^{n} \varphi_{1j,k_t-1}^*\left((k_t-1)T_t\right) x^*\left((k_t-j)T_t\right) \tag{3.14}$$

with the n- real vector $\varphi_{1j,k_t-1}^*$ being defined in (27.a). Eqn. 3.13 may be rewritten through (3.14) for $t = k_t T_t$, $\forall k_t \geq n$ as:

$$x(k_t T_t) = \sum_{j=1}^{n} \varphi_{1j,k_t-1}^*\left((k_t-1)T_t\right) x\left((k_t-j)T_t\right) = c(k_t T_t) := \sum_{j=1}^{n} \left( \varphi_{1j,k_t-1}^*\left((k_t-1)T_t\right) - \varphi_{1j,k_t-1}\left((k_t-1)T_t\right)\right) x^*\left((k_t-j)T_t\right) \tag{3.15}$$

Taking z-transforms in (3.15), one gets for each $t \in \boldsymbol{R}_{0+}$

$$\Delta_{T_t}(z)\, \hat{x}_{T_t}(z) = \hat{q}_{T_t}(z) + \hat{c}_{T_t}(z), \quad |z| > \mu_{T_t} \tag{3.16}$$

with $\mu_{T_t} := \limsup_{k_t \to \infty} \sqrt[k_t]{|x(k_t T_t)|} = \limsup_{t \to \infty} |x(t)|^{T_t / t}$, and $\Delta_{T_t}(z) := \sum_{i=1}^{n} \psi^*_{iT_t} z^{n-i}$ being the characteristic polynomial of $\Psi^*(T_t, 0) = e^{A^* T_t}$, where $\hat{x}_{T_t}(z) = \sum_{k=0}^{\infty} x(kT_t) z^{-k}$ and $\hat{c}_{T_t}(z) = \sum_{k=0}^{\infty} c(kT_t) z^{-k}$ are the Z-transforms of the sequences $\{x(kT_t)\}_0^\infty$ and $\{c(kT_t)\}_0^\infty$, respectively, and

$$\hat{q}_{T_t}(z) := \sum_{k=0}^{n-1} x(kT_t) z^{n-k} + \sum_{j=1}^{n} \sum_{\ell=0}^{n-j-1} \varphi^*_{1j,k_t-1} x(\ell T_t) z^{n-j-\ell} \tag{3.17}$$

for $T_t = t / k_t$ and some $k_t \in \boldsymbol{N}_0$ depends of initial conditions and is zero if $x(kT_t) = 0$, $\forall k \in \overline{n-1} \cup \{0\}$. Note that $\hat{x}_{T_t} : \boldsymbol{C} \to \boldsymbol{C}$ is holomorphic in $|z| > \mu_{T_t}$. Also, $\hat{c}_{T_t} : \boldsymbol{C} \to \boldsymbol{C}$ is holomorphic in $|z| > \nu_{T_t}$ where $\nu_{T_t} := \limsup_{k_t \to \infty} \sqrt[k_t]{|c(k_t T_t)|} = \limsup_{t \to \infty} |c(t)|^{T_t / t} \le \eta_t \mu_t \le \mu_t$.

Take any $\boldsymbol{R}_0 \ni t = k_t T_t$ with $k_t \in \boldsymbol{N}_0$ (properly speaking $\boldsymbol{R}_{0+} \ni T_t := t / k_t < \infty$ for some given injective map $k_t : \boldsymbol{R}_{0+} \to \boldsymbol{N}_0$ $\forall t \in \boldsymbol{R}_{0+}$ satisfying $k(0) = 0$, $k(t) \ne 0$ for $t \ne 0$). Now, choose $\varepsilon_{T_t}$ sufficiently small such that $\mu_{T_t} - 2\varepsilon_{T_t} > \nu_{T_t}$ and all the roots of the characteristic polynomial $\Delta_{T_t}(z)$ belong to the open annulus $\Omega_{T_t} := \{ z \in \boldsymbol{C} : \mu_{T_t} - 2\varepsilon_{T_t} < |z| < \mu_{T_t} + 2\varepsilon_{T_t} \}$. Thus, one gets from (3.16):

$$x(t) = x(k_t T_t) = \frac{1}{2\pi i} \oint_{\gamma_a} z^{k_t - 1} \hat{x}_{T_t}(z) \, dz = \frac{1}{2\pi i} \oint_{\gamma_b} z^{k_t - 1} \Delta^{-1}_{T_t}(z) f(z) \, dz \; ; \; \forall t \in \boldsymbol{R}_{0+} \tag{3.18}$$

where $f(z) = \hat{q}_{T_t}(z) + \hat{c}_{T_t}(z)$, $\gamma_a$ is any positively oriented simple closed curve that lies in the region $|z| > \mu_{T_t}$ and winds around the origin and $\gamma_b : [0, 2\pi] \to \boldsymbol{R}_+$ is the positively oriented circumference centred at the origin defined by $\gamma_b(\omega) := (\mu_{T_t} + \varepsilon_{T_t}) e^{i\omega}$, $\forall \omega \in [0, 2\pi]$. Note that $\hat{c}_{T_t}(z)$ is holomorphic for $|z| > \nu_{T_t}$. Thus, $f(z)$ is holomorphic in $|z| > \nu_{T_t}$ since $\hat{q}_{T_t}(z)$ is a complex-valued polynomial, then an entire complex-valued function. Then, the function $id^{k_t - 1} \Delta^{-1}_{T_t} f : \boldsymbol{C} \to \boldsymbol{C}$ is meromorphic in $\Omega_{T_t}$, $\forall k_t \in \boldsymbol{N}_0$, $\forall t \in \boldsymbol{R}_{0+}$ for $T_t = t / k_t$ with all its poles in $\Lambda(\mu_{T_t})$. Let $\gamma_c : [0, 2\pi] \to \boldsymbol{R}_+$ be the positively oriented circumference centred at the origin defined by $\gamma_c(\omega) := (\mu_{T_t} - \varepsilon_{T_t}) e^{i\omega}$, $\forall \omega \in [0, 2\pi]$ and let $\gamma_d : [0, 2\pi] \to \boldsymbol{R}_+$ be the opposite to $\gamma_c : [0, 2\pi] \to \boldsymbol{R}_+$. Define the cycle $\Gamma := \gamma_b \dot{+} \gamma_c$. By the residue theorem, one gets from (3.18):

$$x(t) \equiv x(k_t T_t) = \frac{1}{2\pi i} \oint_{\gamma_b} z^{k_t - 1} \hat{x}_{T_t}(z) \, dz = \frac{1}{2\pi i} \oint_{\gamma_c} z^{k_t - 1} \Delta^{-1}_{T_t}(z) f(z) \, dz + x^*(t)$$

$$= x^*(t) - \frac{1}{2\pi i}\oint_{\gamma_d} z^{k_t - 1}\Delta_{T_t}^{-1}(z)\, f(z)\, dz \;\; ; \forall t \in \boldsymbol{R}_{0+} \tag{3.19a}$$

since $x^*(t) = \sum_{\lambda \in A_{T_t}} Res\left(id^{k_t - 1}\Delta^{-1} f\right)$ with $A_{T_t} \subset \Lambda_{\mu_t}$ being the set of poles of $id^{k_t - 1}\Delta^{-1} f$ in the annulus $\Omega_{T_t}$, [2] since $\left|\oint_{\gamma_c} z^{k_t - 1}\Delta_{T_t}^{-1}(z)\, f(z)\, dz\right| \le 2\pi K\left(\mu_{T_t} - \varepsilon_{T_t}\right)^{k_t}$ where $K := \max_{|z| = \mu - \varepsilon}\left|\Delta^{-1} f\right|$, [3]. This implies that $x(k_t T_t) = x^*(k_t T_t) + O\left(\left(\mu_{T_t} - \varepsilon_{T_t}\right)^{k_t}\right)$ as $\boldsymbol{R}_0 \ni t = k_t T_t \to \infty$ for any prefixed finite $T_t \in \boldsymbol{R}_+$. Finally, note that $x^*(k_t T_t) \equiv 0$ ; $\forall k_t \in \boldsymbol{N}_0$ and any given $T_t \in \boldsymbol{R}_+$ would imply $\mu_{T_t} = \limsup_{k_t \to \infty} |x(k_t T_t)|^{1/k_t} < \mu_{T_t} - \varepsilon_{T_t} < \mu_{T_t}$ which is a contradiction . Thus, $x^*(t) \equiv x^*(k_t T_t)$ is a nontrivial solution of the limiting equation of (2.1) coincident with the solution of the corresponding limiting discrete equation for any sampling period $T_t \in \boldsymbol{R}_+$. □

REMARK 3.8. The use of a time- dependent sampling period $T_t$ in Theorem 3.7 is a constructive technique to obtain an asymptotic approximation formula for all large time between the solution of (2.1) and its limiting counterpart by using computation residues theory their respective (sampling period–dependent) associate discrete-time difference equations. The sampling period is re-calculated for any t when necessary to obtain a time-dependent nonnegative integer $k_t$ valid to build the identity $t = k_t T_t$ to proceed with the appropriate residues formulas to obtain an asymptotic comparison of the solutions of (2.1) and its associate limiting counterpart. If the sampling period is prefixed to a constant value T then the solution of (2.1) is compared to its limiting solution for large time instants $t = kT$. For $t \ne kT$, Theorem 3.7 is not directly applicable but there is a close version extendable to this case addressed in Corollary 3.9 below where continuity and boundedness arguments for both solutions are used. □

Assume that the differential equation (2.1) and its associated limiting differential equation are both globally asymptotically Lyapunov´ s stable. For this purpose, it suffices from Theorem 2.1 that all the characteristic roots of (2.1), or equivalently, that all the eigenvalues of the matrix A (t) in (2.3) be located in the complex open left half plane Re s < 0. Assume also that the convergence of the coefficient functions of (2.1) to their limits is exponential (Theorem 3.7). Then, $\mu_{T_t} < 1$, $\forall T_t \in \boldsymbol{R}_+$.

Theorem 3.7 becomes a stronger result if (2.1) and the limiting equation are globally Lyapunov´s stable since both solutions converge asymptotically to each other and to zero as it is proven in the subsequent result:

COROLLARY 3.9. Assume that the differential equation (2.1) and its associated limiting differential equation are both globally asymptotically Lyapunov´ s stable. Assume also that the convergence of the coefficient functions of (2.1) to their limits is exponentially fast. Then, $0 < \mu_{T_t} < 1$, $\forall T_t \in \boldsymbol{R}_+$.

*Proof*: Since the differential equation (2.1) and its associated limiting differential equation are both globally asymptotically Lyapunov´s stable, so that $\mu_{T_t} < 1, \forall T_t \in \mathbf{R}_+$ and, furthermore, $x(k_t T_t) = x^*(k_t T_t) + O\left(\left(\mu_{T_t} - \varepsilon_{T_t}\right)^{k_t}\right)$ from Theorem 3.7, it follows that:

$x(k_t T_t) = x^*(k_t T_t) + o\left(\left(\mu_{T_t} - \varepsilon_{T_t}\right)^{k_t}\right)$

as $\mathbf{R}_0 \ni t = k_t T_t \to \infty$ for any prefixed finite $T_t \in \mathbf{R}_+$, $\mathbf{N}_0 \ni k_t \to \infty$, and for some $\varepsilon_{T_t} \in \left(0, \mu_{T_t}\right)$ and some finite $\mu_{T_t} \in \mathbf{R}_+$. Also, since $0 < \mu_{T_t} < 1$ and $0 < \varepsilon_{T_t} < \mu_{T_t} < 1$ then $0 < \mu_{T_t} - \varepsilon_{T_t} < 1$, $\forall T_t \in \mathbf{R}_+$:

$$\lim_{\mathbf{N}_0 \ni k_t \to \infty} x(k_t T_t) = \lim_{\mathbf{N}_0 \ni k_t \to \infty} \left(x^*(k_t T_t) + o\left(\left(\mu_{T_t} - \varepsilon_{T_t}\right)^{k_t}\right)\right) = \lim_{\mathbf{N}_0 \ni k_t \to \infty} x^*(k_t T_t) = 0 \qquad \square$$

A particular version of Corollary 3.9 is the following:

COROLLARY 3.10. Corollary 3.9 holds if all the eigenvalues of the matrix A (t) in (2.3) are located in the complex open left half plane Re s < 0 and convergence of the coefficient functions of (2.1) to their limits is exponentially fast .

Theorem 3.7 may be reformulated for the case of a constant sampling period T as follows:

COROLLARY 3.11. Consider any finite constant sampling period $T \in \mathbf{R}_+$. Then, $x(kT+\tau) = x^*(kT+\tau) + O\left(\left(\mu_T - \varepsilon_T\right)^k\right)$; $\forall \tau \in [0, T)$ as $\mathbf{R}_0 \ni t = kT + \tau \to \infty$ and $\mathbf{N}_0 \ni k \to \infty$, for some $\varepsilon_T \in \left(0, \mu_T\right)$ and some $\mu_T \in \mathbf{R}_+$ such that $\Lambda(\mu_T) := \left\{ z \in \mathbf{C} : \Delta_{T_t}(z) = 0, |z| = \mu_T \right\} \neq \emptyset$. Then, $x(k_t T_t) = x^*(k_t T_t) + O\left(\left(\mu_{T_t} - \varepsilon_{T_t}\right)^{k_t}\right)$ as $\mathbf{R}_0 \ni t = k_t T_t \to \infty$ for any prefixed finite $T_t \in \mathbf{R}_+$, $\mathbf{N}_0 \ni k_t \to \infty$, and for some $\varepsilon_{T_t} \in \left(0, \mu_{T_t}\right)$ and some $\mu_{T_t} \in \mathbf{R}_+$ such that

*Proof*: Note from Theorem 3.7 and using the property

$$f_1(t) = O(f_2(t)) \Leftrightarrow |f_1(t)| \le K_1 |f_2(t)| + K_2$$

for some finite $K_i \in \mathbf{R}_{0+}$ (i=1,2) of the Landau´s big – O, that for $\mathbf{R}_0 \ni t = kT + \tau \to \infty$:

$$\begin{aligned}
\left|x(kT+\tau) - x^*(kT+\tau)\right| &= \left|\left(x(kT) - x^*(kT)\right) + \left(x(kT+\tau) - x(kT)\right) + \left(x^*(kT) - x^*(kT+\tau)\right)\right| \\
&\le \left|x(kT) - x^*(kT)\right| + \left|x(kT+\tau) - x(kT)\right| + \left|x^*(kT) - x^*(kT+\tau)\right| \\
&\le O\left(\left(\mu_T - \varepsilon_T\right)^k\right) + K_1(\tau) + K_1(\tau) \le O\left(\left(\mu_T - \varepsilon_T\right)^k\right) + \overline{K}_1 + \overline{K}_2 \\
&\le O\left(\left(\mu_T - \varepsilon_T\right)^k\right) + \overline{K} = O\left(\left(\mu_T - \varepsilon_T\right)^k\right)
\end{aligned}$$

where:

$\overline{K} := \overline{K}_1 + \overline{K}_2$ , with $\overline{K}_1 := \max_{\tau \in [0, T)} K_1(\tau)$ and $\overline{K}_2 := \max_{\tau \in [0, T)} K_2(\tau)$

are real finite nonnegative constants depending on T, but independent of $\tau$ , since the solution of the linear differential equation (2.1) and that of its limiting equation are both continuously differentiable then uniformly bounded within any finite real interval provided that they are bounded at any point of such an interval . □

ACKNOWLEGEMENTS

The author is very grateful to the Spanish Ministry of Education by its partial financial support of this work through Project DPI2006-0714. He is also grateful to the Basque Government by its partial financial support through Grants SAIOTEK S-PE07UN04 and GIC 07143-IT-269-07.